\documentclass{article}
\usepackage[affil-it]{authblk}
\usepackage{amsthm, graphicx, latexsym}  
\usepackage[utopia]{mathdesign}      
 \usepackage[mathscr]{eucal}
\usepackage{amsmath, yhmath}  
\usepackage{array, longtable, sidecap, wrapfig, subfig, color}
\usepackage[normalem]{ulem}
\usepackage{cancel}
\usepackage{tikz}
\usepackage{stmaryrd}  

\usepackage[colorlinks=true, pdfstartview=FitV, linkcolor=blue, citecolor=blue, urlcolor=cerulean]{hyperref}
\DeclareMathOperator{\fix}{fix}
\DeclareMathOperator{\dom}{dom}
\DeclareMathOperator{\sgn}{sgn} 

\usepackage{url}
\usepackage{fancyhdr}
\begin{document}

\newtheorem{thm}{Theorem}
\newtheorem*{thm*}{Theorem}
\newtheorem{cor}[thm]{Corollary}
\newtheorem{lemma}[thm]{Lemma}
\newtheorem{lemmadef}[thm]{Lemma/Definition}
\newtheorem{proposition}[thm]{Proposition}
\newtheorem{prop}[thm]{Proposition}
\newtheorem{conj}[thm]{Conjecture}
\newtheorem{example}[thm]{Example}

\newenvironment{definition}[1][Definition]{\begin{trivlist}
\item[\hskip \labelsep {\bfseries #1}]}{\end{trivlist}}
\newenvironment{remark}[1][Remark]{\begin{trivlist}
\item[\hskip \labelsep {\bfseries #1}]}{\end{trivlist}}
\newenvironment{remarks}[1][Remarks]{\begin{trivlist}
\item[\hskip \labelsep {\bfseries #1}]}{\end{trivlist}}

\newenvironment{solution}{\paragraph{Solution:}}{\hfill QEF.}

\theoremstyle{definition}
\newtheorem{defn}[thm]{Definition}
\newtheorem{defns}[thm]{Definitions}
\newtheorem{rem}[thm]{Remark}
\newtheorem{Rem}[thm]{Remark}
\newtheorem{Rems}[thm]{Remarks}
\newtheorem{Exercise}[thm]{Exercise}
\newtheorem{Examples}[thm]{Examples}
\newtheorem{Problem}[thm]{Problem}
\newtheorem{Open questions}[thm]{Open questions}
\newtheorem{Open question}[thm]{Open question}
\newtheorem{Open problems}[thm]{Open problems}
\newtheorem{Open problem}[thm]{Open problem}
\newtheorem*{Acknowledgements}{Acknowledgements}  
\newtheorem*{Organization}{Organization of the paper}  
\newtheorem*{Asymptotic Notation}{Asymptotic Notation}
\newtheorem*{Generic Complexity}{Generic Complexity}

\newcommand{\leftexp}[2]{{\vphantom{#2}}^{#1}\!{#2}}
\newcommand{\beq}{\begin{eqnarray*}}
\newcommand{\eeq}{\end{eqnarray*}}
\newcommand{\ol}{\overline}
\newcommand{\Supp}{\textrm{Supp}}
\newcommand{\Span}{\textrm{Span}}
\newcommand{\Gr}{\textbf{Gr} }
\newcommand{\GL}{\textrm{GL}}
\newcommand{\SL}{\textrm{SL}}
\newcommand{\Matr}{\textbf{Matr}}
\newcommand{\mat}[9]{\left(\begin{array}{ccc}#1&#2&#3\\#4&#5&#6\\#7&#8&#9\end{array}\right)}
\newcommand{\rowm}[3]{\left(\begin{array}{c}#1\\#2\\#3\end{array}\right)}
\newcommand{\rk}{\textrm{rk}~}
\newcommand{\im}{\textrm{im}~}
\newcommand{\Div}{\textrm{div}~}
\newcommand{\z}{\overline}
\newcommand{\Par}[2]{\frac{\partial#1}{\partial#2}}
\newcommand{\hf}{\frac{1}{2}}
\newcommand{\lf}{\left(}
\newcommand{\rt}{\right)}
\newcommand{\re}{\textrm{Re}}
\newcommand{\img}{\textrm{Im}}
\newcommand{\limn}{{\lim\atop{n\to\infty}}}
\newcommand{\Res}{\textrm{Res}}
\newcommand{\eps}{\epsilon}
\newcommand{\tp}{\frac{1}{2\pi i}}
\newcommand{\of}{\circ}
\newcommand{\ub}{\underbrace}
\newcommand{\K}{\tilde{k}}
\newcommand{\floor}[1]{{\lfloor#1\rfloor}}
\newcommand{\Gal}{\textrm{Gal}}
\newcommand{\Aut}{\textrm{Aut}}
\newcommand{\out}[1]{{\textrm{Out}(F_#1)}}
\newcommand{\pr}[2]{\langle#1,#2\rangle}
\newcommand{\Tr}{\textrm{Tr}}
\newcommand{\id}{\mathfrak}
\newcommand{\rad}{\mathfrak{R}}
\newcommand{\Spec}{\mathrm{Spec}}
\newcommand{\tensor}{\otimes}
\newcommand{\Ann}{\textrm{Ann}}
\renewcommand{\phi}{\varphi}
\newcommand{\Ker}{\textrm{Ker}~}
\newcommand{\Img}{\textrm{Im}~}
\newcommand{\Z}{\mathbb{Z}}
\newcommand{\N}{\mathbb{N}_{\geq0}}
\newcommand{\trans}[1]{{\leftexp{T}{#1}}}
\def\onto{{\kern3pt\to\kern-8pt\to\kern3pt}}
\newcommand{\IO}{\textrm{IO}}
\newcommand{\IA}{\textrm{IA}}
\newcommand{\Stab}{\textrm{Stab}}
\newcommand{\ssim}[1]{\stackrel{(#1)}{\sim}}
\newcommand{\Lk}[1]{\textrm{Lk}_<(#1)}
\newcommand{\set}[1]{\left\{#1\right\}}
\newcommand{\ssm}{\smallsetminus}
\newcommand{\abs}[1]{\left|#1\right|}
\newcommand{\Dist}{\textup{Dist}}
\newcommand{\sem}[1]{\llbracket#1\rrbracket}

\newcommand{\wh}{\widehat}
\newcommand{\wt}{\widetilde}
\newcommand{\wb}{\overline}
\newcommand{\ihat}{\hat{\imath}}
\newcommand{\fhat}{\hat{f}}
\newcommand{\ibar}{\overline{\imath}}
\newcommand{\fbar}{\overline{f}}
\newcommand{\partialto}{\rightharpoonup}
\newcommand{\bigos}{\widetilde{\mathcal{O}}} 
\newcommand{\bigo}{\mathcal{O}}
\newcommand{\bigom}{\Omega}
\newcommand{\bigt}{\Theta}
\newcommand{\emptystring}{\epsilon}
\newcommand{\evale}{\varepsilon}
\newcommand{\redu}{\widetilde}

\def\ms{\medskip}
\newcommand{\bs}{\bigskip}

\newcommand{\tc}[2]{\textcolor{#1}{#2}}
\definecolor{cerulean}{rgb}{0,.48,.65} \newcommand{\cerulean}[1]{\tc{cerulean}{#1}}
\definecolor{magenta}{rgb}{.5,0,.5} \newcommand{\magenta}[1]{\tc{magenta}{#1}}
\definecolor{dred}{rgb}{.5,0,0} \newcommand{\dred}[1]{\tc{dred}{#1}}
\definecolor{green}{rgb}{0,.5,0} \newcommand{\green}[1]{\tc{green}{#1}}
\definecolor{blue}{rgb}{0,0,0.5} \newcommand{\blue}[1]{\tc{blue}{#1}}
\definecolor{black}{rgb}{0,0,0} \newcommand{\black}[1]{\tc{black}{#1}}
\definecolor{dgreen}{rgb}{0,.3,0} \newcommand{\dgreen}[1]{\tc{dgreen}{#1}}
\definecolor{vdred}{rgb}{.3,0,0} \newcommand{\vdred}[1]{\tc{vdred}{#1}}
\definecolor{red}{rgb}{1,0,0} \newcommand{\red}[1]{\tc{red}{#1}}
\definecolor{salmon}{rgb}{0.98,0.50,0.45} \newcommand{\salmon}[1]{\tc{salmon}{#1}}
\definecolor{gray}{rgb}{.5,.5,.5} \newcommand{\gray}[1]{\tc{gray}{#1}}
\definecolor{seagreen}{rgb}{0.13,0.70,0.67} \newcommand{\seagreen}[1]{\tc{seagreen}{#1}}
\definecolor{chartreuse}{rgb}{0.40,0.80,0.00}\newcommand{\chartreuse}[1]{\textcolor{chartreuse}{#1}}
\definecolor{cornflower}{rgb}{0.39,0.58,0.93} \newcommand{\cornflower}[1]{\textcolor{cornflower}{#1}}
\definecolor{gold}{rgb}{0.80,0.68,0.00}\newcommand{\gold}[1]{\textcolor{gold}{#1}}

\setlength{\parindent}{0pt}
\setlength{\parskip}{7pt}

\title{The Conjugacy Problem for Higman's Group}
\author{Owen Baker}
 
\date \today

\maketitle

\begin{abstract}
\noindent In 1951, Higman\cite{Higman} constructed a remarkable group
$$H=\left\langle a,b,c,d \, \left| \,  b^a = b^2, c^b = c^2, d^c = d^2, a^d = a^2 \right. \right\rangle$$
and used it to produce the first examples of infinite simple groups.
By studying fixed points of certain finite state transducers, we show the conjugacy problem in $H$ is decidable (for all inputs).

Diekert, Laun \& Ushakov\cite{DLU2} have recently shown the word problem in $H$ is solvable in polynomial time,
using the \emph{power circuit} technology of Myasnikov, Ushakov \& Won\cite{MUW}.
Building on this work, we show in a strongly generic setting that the conjugacy problem has a polynomial time solution.

\ms \noindent   \textbf{2010 Mathematics Subject Classification:  20F10, 68Q70 }  
\end{abstract}

In 1911--1912, Max Dehn introduced the word problem and the conjugacy problem\cite{Dehn1}\cite{Dehn2}.
Fix a finitely generated group $G$.
The \emph{word problem} asks, given words $x,y$ in the generators and their inverses, whether $x$ and $y$ represent the same element of $G$.
In this case, we write $x=_Gy$.
The \emph{conjugacy problem} asks, given words $x,y$ in the generators and their inverses, whether $x$ and $y$ represent conjugate elements of $G$.  That is, does there exist $z\in G$ so that $x^z:=z^{-1}xz=_Gy$?  In this case, we write $x\sim_Gy$ and call $z$ (or a word representing it) a \emph{conjugator}.
When convenient, we conflate a word with the group element it represents.

The group
\begin{align*}
H   \ &= \ \left\langle a,b,c,d \, \left| \,  b^a = b^2, c^b = c^2, d^c = d^2, a^d = a^2 \right. \right\rangle\\
&=A\ast_CB=\left\langle a,b,c \, \left| \,  b^a
= b^2, c^b = c^2 \right. \right\rangle\ast_{\langle a,c\rangle}
\left\langle c,d,a \, \left| \,  d^c = d^2, a^d = a^2 \right. \right\rangle\\
&=D\ast_FE=\left\langle b,c,d \, \left| \,  c^b = c^2, d^c = d^2 \right. \right\rangle\ast_{\langle b,d\rangle}
\left\langle d,a,b \, \left| \,   a^d = a^2, b^a = b^2 \right. \right\rangle
\end{align*}
was constructed by Graham~Higman in his celebrated 1951 paper \cite{Higman} as an example of a group which has no finite quotients.
The groups $H/N$, with $N$ any maximal proper normal subgroup, were consequently the first examples of finitely generated infinite simple groups.  Here $A\cong B\cong D\cong E$ is an HNN--extension of the Baumslag--Solitar group
$BS(1,2)=\left\langle b,c \, \left| \,  c^b = c^2 \right. \right\rangle$.
$H$ is the amalgamated free product (in two different ways!) of two copies of this group along a rank 2 free group $C\cong F$.
We prove:

\begin{thm} \label{conjugacy in H}
The conjugacy problem in $H$ is decidable.
\end{thm}

Our approach is based on the method that Diekert, Myasnikov, and Wei{\ss} applied to the conjugacy problem for the Baumslag group $G_{1,2}$ in \cite{DMW2, DMW1}.  The group $G_{1,2}$ is also known as the Baumslag--Gersten group, and its conjugacy problem was previously solved by Beese in \cite{Beese}, a German Diploma thesis.  However, there is a complication.
\cite{DMW2, DMW1} and \cite{Beese} ultimately reduce the problem ``is $x\sim_Gy$?'' to testing whether $x^z=_Gy$ for a single candidate 
 conjugator $z$, depending on $x$ and $y$.  Diekert et al.\ call this the ``key'' to their approach.

For Higman's group, we find cases where there is an infinite family of candidates to check.
To deal with this, we build from $x,y$ a \emph{finite state automaton} which decides if a given input word serves as a conjugator (section \ref{Section:Equation}).
One then algorithmically checks if there is a path from the start node to any accept node.
If so, the path's label gives a conjugator $z$.  If not, $x\nsim_Hy$.
Thus we are able to algorithmically solve the conjugacy problem in $H$.

The word problem for Higman's group $H$ is solvable in polynomial time.
Specifically, Diekert, Laun \& Ushakov gave a $\bigo(n^6)$ algorithm in \cite{DLU2} based on the \emph{power circuit} technology of Myasnikov, Ushakov, \& Won\cite{MUW}.
(The journal version \cite{DLU2} improves on an earlier version \cite{DLU1}, using amortized analysis to remove logarithmic factors.
Clarification of the amortized analysis is provided by Laun\cite[section 2.4.2]{Laun}.)
Prior to this result, $H$ had been a candidate for a group with hard but decidable word problem.
Indeed, the Dehn function of $H$ is non-elementary, with a lower bound involving the tower of exponents function.
Decidability of the word problem follows, without \cite{DLU2}, from the fact that $H$ is iteratively built from $\mathbb{Z}$ by HNN--extensions and amalgamated free products.

Borovik, Myasnikov \& Remeslennikov show in \cite{BMR} that even in amalgamated free products with undecidable conjugacy problem,
there can be efficient (e.g.\ polynomial time) solutions on \emph{generic} inputs.
Diekert, Myasnikov, and Wei{\ss} give a \emph{strongly generic} quartic time algorithm, using power circuits,
 for the conjugacy problem of $G_{1,2}$ in \cite{DMW2, DMW1, DMW3}.  Roughly speaking, this means the proportion of inputs for which
the algorithm fails to yield an answer decays exponentially with input size.
We show:
\begin{thm}\label{strongly generic conjugacy in H}
Let $\Sigma=\{a^\pm,b^\pm,c^\pm,d^\pm\}$.
There is a strongly generic algorithm that decides in time $\bigo(n^7)$ on input words $x,y\in\Sigma^*$ with total length $n$
whether $x\sim_Hy$.
The algorithm is also strongly generic in time $\bigo(n^7)$ on freely reduced inputs and on cyclically reduced inputs.
\end{thm}

\begin{Asymptotic Notation}
We use standard Big--O notation $\bigo(f)$, as well as $\bigt(f)$ and $\bigom(f)$.
\end{Asymptotic Notation}

\begin{Generic Complexity}
The notion of (strong) generic complexity was introduced in \cite{GenericCaseComplexity}.
Consider an algorithm $\mathcal{A}$ taking inputs from a domain $D$.
For us, $D$ is either $\Sigma^*$ (the set of finite words in $\Sigma$); the subset of freely reduced or of cyclically reduced words;
or the set of pairs of such words.
In each case, there is a natural partition $D=\amalg_{n\in\mathbb{N}}D^{(n)}$ into finite sets $D^{(n)}$ of ``size $n$ inputs''.
For example, $D^{(n)}$ may be the set $\Sigma^n\subseteq\Sigma^*$ of length $n$ words,
or the set $\bigcup_{i=0}^n{\Sigma^i\times\Sigma^{n-i}}\subseteq\Sigma^*\times\Sigma^*$ of pairs of total length $n$.

A set $I\subseteq D$ is called \emph{generic} if $|I\cap D^{(n)}|/|D^{(n)}|\to1$ as $n\to\infty$, and \emph{strongly generic} if
$$\frac{|I\cap D^{(n)}|}{|D^{(n)}|}=1-e^{-\bigom(n)}.$$
The complement of a (strongly) generic set is \emph{(strongly) negligible}.

The algorithm $\mathcal{A}$ \emph{runs in (strongly) generic time $\bigo(f)$} if there is a (strongly) generic set $I\subseteq D$
so that $\mathcal{A}$ takes at most $\bigo(f(n))$ steps on each element of $|I\cap D^{(n)}|$.
Though $\mathcal{A}$ need not terminate on input outside $I$, it must never halt with incorrect answers.

In an abuse of language, we say ``a (strongly) generic $w\in D$ satisfies predicate $Q(w)$'' to mean $\{w:Q(w)\}$ is (strongly) generic in $D$.
\emph{Of course, an individual element cannot be classified as ``generic'' or not.}
\end{Generic Complexity}

\begin{Organization}
Section \ref{Section:Warmup} reviews background about the conjugacy problem for HNN--extensions and free products,
and as a warm-up applies the techniques of \cite{DMW1}\cite{DMW2} to solve the conjugacy problem for $A\cong B\cong D\cong E$.
Section \ref{Section:Automata} reviews background about finite state automata, and develops a criterion under which a fixed point set forms a regular language.
Section \ref{Section:Equation} solves the conjugacy problem for Higman's group $H$.
Section \ref{Section:Runtime} shows that the algorithm runs in polynomial time on a large class of inputs, by applying the \emph{power circuit} data compression techniques of \cite{DLU2}.
\end{Organization}

\begin{Acknowledgements}
The author thanks Timothy Riley for recommending this problem and explaining his elegant geometric intuition underlying conjugacy in Higman's group.
\end{Acknowledgements}

\section{Preliminaries and a warm--up conjugacy problem}\label{Section:Warmup}
We assume the reader is familiar with HNN--extensions and amalgamated free products; see the classic textbooks\cite{LS}\cite{MKS}.
In this section, we establish terminology and review some of the standard definitions and facts.
As a warm-up, we prove the following proposition, illustrating the \cite[Section 5]{DMW1} techniques.

\begin{prop} \label{conjugacy in A}
The conjugacy problem is solvable in $A\cong B\cong D\cong E$ where, we recall,
$$A=\left\langle a,b,c \, \left| \,  b^a = b^2, c^b = c^2 \right. \right\rangle.$$
\end{prop}

Here $A$ is an HNN--extension\footnote{Several authors investigate solvability of the conjugacy problem generally in groups splitting over
     cyclic edge groups (\cite{Hurwitz}\cite{HoradamFarr} for HNN--extensions, \cite{Lipschutz} for amalgamated free products, 
		 and \cite{Horadam} for general graphs of groups).
     These results would here require centrality of $b$ in $K$, or else a \emph{semi-criticality} condition implying $b\nsim_Ab^2$.
     So they do not seem to help prove Proposition \ref{conjugacy in A}.
	}
over the Baumslag--Solitar group $K=\left\langle b,c\,\left|\, c^b=c^2\right.\right\rangle$
with stable letter $a$ and associated subgroups $\langle b\rangle$ and $\langle b^2\rangle$.
We have $K\cong\mathbb{Z}[1/2]\rtimes\mathbb{Z}$ via the correspondence $c=(1,0)$ and $b=(0,1)$.
Elements of $\mathbb{Z}[1/2]\rtimes\mathbb{Z}$ are pairs $(r,s)$ where $r$ is a dyadic rational and $s$ is an integer.
The group structure is given by $(r,s)(r',s')=(r+2^{-s}r',s+s')$.
\textbf{Caution: }We use slightly different notation from \cite{DMW2, DMW1}, since for us $c^b$ means $b^{-1}cb$ rather than $bcb^{-1}$.

Each element of the HNN--extension $A$ is represented (non-uniquely) by a word $w$ in the (infinite) alphabet $\{a^{+1},a^{-1}\}\amalg\mathbb{Z}[1/2]\rtimes\mathbb{Z}$.
\begin{definition}
Deleting a substring $a^{\pm1}a^{\mp1}$ from $w$ or replacing $(r,s)(r',s')$ with its product in $\mathbb{Z}[1/2]\rtimes\mathbb{Z}$ is called \emph{reduction}.
Replacing a substring of the form $a^{-1}(0,s)a^{+1}$ by $(0,2s)$ is called \emph{Britton--reduction}, as is replacing
$a^{+1}(0,2s)a^{-1}$ by $(0,s)$.
If no reduction or Britton--reduction is possible, we say $w$ is \emph{Britton--reduced}.
We say $w$ is \emph{cyclically Britton--reduced} if $ww$ is \emph{Britton--reduced}.
(If $w=(r,s)$ is a single letter, we also consider $w$ to be cyclically Britton--reduced by fiat.)
\end{definition}

Since Britton--reduction is effective, decidability of the word problem in $A$ follows from:
\begin{lemma}\label{Britton HNN}(Britton's Lemma \cite[IV.2.1]{LS})
A non-empty Britton--reduced word $w$ represents $1\in A$ if and only if $w=(0,0)$ is the identity letter.
More generally, $w$ represents an element of the subgroup $K$ if and only if $w=(r,s)$ is a single letter.
\end{lemma}

To decide if two words $u,v$ represent conjugate elements of $A$, it is convenient to first replace them by cyclically reduced words
$\widehat{u},\widehat{v}$ for which $u\sim_A\widehat{u}$ and $v\sim_A\widehat{v}$.
Of course, any word is conjugate to the Britton--reductions of its cyclic permutations.
The following observation is useful here:
\begin{lemma}\label{cyclic reduction HNN}(\cite[Remark 3, p.973--974]{DMW1})
Given a word $u$ in $\{a^{+1},a^{-1}\}\amalg\mathbb{Z}[1/2]\rtimes\mathbb{Z}$, there is always a cyclic permutation whose Britton--reduction $\widehat{u}$ is cyclically Britton--reduced.
\end{lemma}

Our main tool for proving Proposition \ref{conjugacy in A} is:
\begin{lemma}\label{Collins HNN}(Collins' Lemma\cite[IV.2.5]{LS})
Suppose $u$ is cyclically Britton--reduced, and is not a single letter $(r,s)$.
Then every cyclically Britton--reduced conjugate of $u$ is equal in $A$ to the result of first cyclically permuting the letters of $u$ appropriately and then conjugating the result by an appropriate element of $\langle b\rangle$.
\end{lemma}

\begin{lemma}\label{A recurrence}
Let $x=a^{\epsilon_1}g_1a^{\epsilon_2}g_2\cdots a^{\epsilon_n}g_n$ and $y=a^{\epsilon_1}h_1a^{\epsilon_2}h_2\cdots a^{\epsilon_n}h_n$
be cyclically Britton--reduced, with $\epsilon_i=\pm1$ and $g_i,h_i\in K$.
Then $x,y$ are conjugate in $A$ by an element of $\langle b\rangle$ if and only if the equation
\begin{equation}\label{A equation}
h_{i+1}b_{i+1}=_K\varphi^{\epsilon_{i+1}}(b_i)g_{i+1},\qquad 0\leq i<n
\end{equation}
has a solution sequence $(b_j)_{j=0}^n$ satisfying the boundary condition
\begin{equation}\label{A boundary}
b_0=b_n.
\end{equation}
\end{lemma}
\begin{proof}
Define a homomorphism $\varphi:\langle b\rangle\to\langle b^2\rangle$ so that $a^{-1}ba=\varphi(b)=b^2$.

If $b_0xb_0^{-1}=y$, then the following calculation shows equation (\ref{A equation}) has a recursive solution satisfying (\ref{A boundary}).
\begin{align*}
a^{\epsilon_1}h_1a^{\epsilon_2}h_2\cdots a^{\epsilon_n}h_n&=b_0\left(a^{\epsilon_1}g_1a^{\epsilon_2}g_2\cdots a^{\epsilon_n}g_n\right)b_0^{-1}\\
&=a^{\epsilon_1}\left[\varphi^{\epsilon_1}(b_0)g_1\right]a^{\epsilon_2}g_2\cdots a^{\epsilon_n}g_nb_0^{-1}\\
&=a^{\epsilon_1}\left[h_1b_1\right]a^{\epsilon_2}g_2\cdots a^{\epsilon_n}g_nb_0^{-1}\\
&=a^{\epsilon_1}h_1a^{\epsilon_2}\left[\varphi^{\epsilon_2}(b_1)g_2\right]\cdots a^{\epsilon_n}g_nb_0^{-1}\\
&=a^{\epsilon_1}h_1a^{\epsilon_2}\left[h_2b_2\right]\cdots a^{\epsilon_n}g_nb_0^{-1}\\
&=\cdots\\
&=a^{\epsilon_1}h_1a^{\epsilon_2}h_2\cdots a^{\epsilon_n}\left[h_nb_n\right]b_0^{-1}\\
&=a^{\epsilon_1}h_1a^{\epsilon_2}h_2\cdots a^{\epsilon_n}h_n\left(b_nb_0^{-1}\right)
\end{align*}
The converse follows by the same computation.
\end{proof}

\begin{lemma}\label{A recurrence 2}
Let $x,y$ be as in Lemma \ref{A recurrence}.
Fix some specific $I$, $0\leq I<n$.
Let $h_{I+1}=(r,s)$ and $g_{I+1}=(r',s')$.
If $r\neq0$ or $r'\neq0$, there is at most one sequence satisfying equation (\ref{A equation}) -- without necessarily satisfying (\ref{A boundary}).
In this case, we can calculate $b_0,\ldots,b_n$ effectively (or conclude no such sequence exists) -- and in particular decide whether
$b_0xb_0^{-1}=y$.
\end{lemma}
\begin{proof}
Write $b_{I+1}=b^n$ and $\varphi^{\epsilon_{I+1}}(b_I)=b^{n'}$.  Then equation (\ref{A equation}) becomes:
$$(r,s)(0,n)=(0,n')(r',s')\iff r=2^{-n'}r'\textrm{ and } s+n=s'+n'.$$
Now, $r=2^{-n'}r'$ has at most one solution $n'$, unless $r=r'=0$.  Once $n'$ is known, $s+n=s'+n'$ can be solved for $n$.
Once $b_{I+1}=b^n$ is known, equation (\ref{A equation}) can be solved recursively for all other $b_i$.
Finally, by Lemma \ref{A recurrence}, we have $b_0xb_0^{-1}=y$ if and only if $b_0=b_n$.
\end{proof}

We now solve the conjugacy problem in $A$:

\begin{proof}[Proof of Proposition \ref{conjugacy in A}]
Given words $x,y$ in $\{a^\pm,b^\pm,c^\pm\}$, we must decide whether $x\sim_Ay$.
Replace $x,y$ by words in the alphabet $\{a^+,a^-\}\amalg\mathbb{Z}[1/2]\rtimes\mathbb{Z}$
and then by cyclically Britton--reduced words as in Lemma \ref{cyclic reduction HNN}. 
We can check if $x,y\in K$ by Lemma \ref{Britton HNN}.
We first handle the case $x\notin K$.

Collins' Lemma reduces us to the case where $x=a^{\epsilon_1}g_1a^{\epsilon_2}g_2\cdots a^{\epsilon_n}g_n$ and 
                                            $y=a^{\epsilon_1}h_1a^{\epsilon_2}h_2\cdots a^{\epsilon_n}h_n$
for some $\epsilon_i=\pm1$ and $g_i,h_i\in K$; we wish to decide whether $b_0xb_0^{-1}=y$ for some $b_0\in\langle b\rangle$.
Lemmas \ref{A recurrence} and \ref{A recurrence 2} achieve this unless all $g_i,h_i\in\langle b\rangle$.

So suppose all $g_i,h_i\in\langle b\rangle$.  View $x,y$ as words in $\{a^\pm,b^\pm\}$.
The quotient $A/\ll c\gg\cong\left\langle a,b \, \left| \,  b^a = b^2 \right. \right\rangle=:BS$.
So if $x\sim_Ay$ then $x\sim_{BS}y$.
Conversely, if $x\sim_{BS}y$ then $x\sim_Ay$ because 
$$BS=\left\langle a,b \, \left| \,  b^a = b^2 \right. \right\rangle
\leq\left\langle a,b \, \left| \,  b^a = b^2 \right. \right\rangle
       \ast_{\langle b\rangle} {\left\langle b,c \, | \,  c^b = c^2 \right\rangle} =A.$$
Thus we can decide whether $x\sim_Ay$ by an appeal to the conjugacy problem in the Baumslag--Solitar group $BS$,
 solved in \cite[Theorem 2]{DMW1}.
This completes the proof in the case $x\notin K$.

There remains the case $x,y\in K$.
Let $x=:(r,m)$ and $y=:(s,q)$.
We must decide if they are conjugate in $A$.
Let $g\in K$ and $\alpha=a^{\pm1}$.
Then $\alpha g(r,m)g^{-1}\alpha^{-1}$ Britton--reduces only if $(r,m)\sim_{K}b^m$.
As long as $(r,m)\nsim_{K}b^m$, we therefore have $(r,m)\sim_A(s,q)\iff(r,m)\sim_{K}(s,q)$.
This latter condition is decidable:
   the conjugacy problem for the Baumslag--Solitar group $K$ is solved in \cite[Theorem 2]{DMW1}, as mentioned above.

We are left with the question: When is $(0,m)\sim_A(0,q)$?
Answer: if $q=m\cdot 2^i$ then $b^m=a^{i}b^qa^{-i}$ so that $b^m\sim_A b^q$.
Conversely, if $b^m\sim_A b^q$ then $b^m$ is conjugate to $b^q$ in $A/\ll c\gg=\langle a, b | b^a=b^2\rangle=BS$, so $q=m\cdot 2^i$
for some $i\in\mathbb{Z}$.
\end{proof}

We close this section with analogous definitions and tools for the amalgamated free product $H=A\ast_CB$.

\begin{definition}
Consider a word $w$ in the (infinite) alphabet $A\amalg B$ representing an element of $H=A\ast_CB$.
Here, \emph{reduction} means combining adjacent letters from the same factor $A$ or $B$ by multiplication therein.
\emph{Britton reduction} means replacing a letter in $C\subset A$ by the corresponding element $C\subset B$, or vice versa, and then reducing.
If neither reduction nor Britton reduction is possible, we say $w$ is \emph{Britton--reduced}.
We say $w$ is \emph{cyclically Britton--reduced} if $ww$ is \emph{Britton--reduced}.
(By fiat, any $w\in A\amalg B$ is cyclically Britton--reduced.)
\end{definition}

\begin{lemma}\label{Britton amalgams}(Britton's Lemma \cite[IV.2.6]{LS})
A non-empty Britton--reduced word $w$ represents an element of $A\cup B$ if and only if it is a single letter.
In particular, $w=_H1$ if and only if $w$ is the letter $1\in A$ or $1\in B$.
\end{lemma}

\begin{lemma}\label{Britton effective}(\cite[Lemma 20]{DLU2})
There is a procedure which, given $w\in\{a^\pm,b^\pm,c^\pm\}^*$,
finds $w'\in\{a^\pm,c^\pm\}^*$ such that $w=_Aw'$, or declares no such $w'$ exists.
Therefore, (cyclic/regular) Britton reduction in $H$ is effective.
\end{lemma}

The same proof as Lemma \ref{cyclic reduction HNN} gives:
\begin{lemma}\label{cyclic reduction amalgams}
Given a word $u$ in $A\amalg B$, there is always a cyclic permutation whose Britton--reduction $\widehat{u}$ is cyclically Britton--reduced.
\end{lemma}

The analogue of Collins' Lemma is:
\begin{lemma}\label{conjugation for amalgams} (\cite[IV.2.8]{LS}; \cite[4.6]{MKS}) 
Let $H=A\ast_CB$ and let $u=a_1b_2\cdots a_{2n-1}b_{2n}\in H$, $n\geq1$, be cyclically Britton--reduced.
Every cyclically Britton--reduced conjugate of $u$ is equal in $H$ to the result of first cyclically permuting $a_1b_2\cdots a_{2n-1}b_{2n}$ appropriately and then conjugating by an appropriate element of $C$.
\end{lemma}

In section \ref{Section:Equation}, we will solve the analogue for $H$ to equation (\ref{A equation}):
Given $a_1,a_1'\in A$, find all $(\gamma_1,\gamma_2)\in C^2$ so that $\gamma_1a_1=a_1'\gamma_2$.
We will use \emph{automata}, the subject of section \ref{Section:Automata}, to assist with the harder cases.

\section{Automata}\label{Section:Automata}
For a background on finite state automata in group theory, see the textbook \cite{ECHLPT}.
Recall:

\begin{definition}
\vbox{%
A \emph{(partial deterministic) automaton} $M$ consists of the data $(S,\Sigma,\mu,s_0,Y)$, where
\begin{itemize}
\item $S$ is a finite \emph{state set};
\item $\Sigma$ is a finite \emph{alphabet};
\item $s_0\in S$ is the \emph{initial state};
\item $\mu:S\times\Sigma\partialto S$, a partial function, is the \emph{transition function}; and
\item $Y\subseteq S$ is the set of \emph{accept states}.
\end{itemize}}
A \emph{language} over $\Sigma$ is a subset of the set $\Sigma^*$ of finite strings.
The data $(S,\mu)$ can be interpreted as a directed graph with vertex set $S$ and arrows labelled by elements of $\Sigma$.
The transition function induces $S\times\Sigma^*\partialto S$, also denoted $\mu$, by the recursion $\mu(s,uv)=\mu(\mu(s,u),v)$ where $s\in S$ and $u,v\in\Sigma^*$.
For the empty string $\emptystring\in\Sigma^*$, we take $\mu(s,\emptystring)=s$.
A word $w\in\Sigma^*$ is \emph{accepted} by $M$ if $\mu(s_0,w)\in Y$.
The set of words accepted by $M$ is denoted $L(M)$.
Thus $w\in L(M)$ if and only if the path starting from $s_0$ and labelled by $w$ terminates in $Y$.
A language is \emph{regular} if it is of the form $L(M)$ for some partial deterministic automaton $M$.
\end{definition}

\begin{lemma}\label{emptiness problem}(\cite[Theorem 7]{RabinScott})
There is an effective procedure deciding, given $M$, whether $L(M)$ is empty.
\end{lemma}

\begin{proof}Check if there exists a path from $s_0$ to any element of $Y$.\end{proof}

We wish to solve the conjugacy problem.  To tell if $x\sim y$, we will want to decide if $z\mapsto y^{-1}zx$ has any fixed points.
To do this, we will build an automaton accepting the fixed point set and then apply Lemma \ref{emptiness problem}.
To study graphs of functions, we need an automaton analogue for \emph{2-variable languages}, subsets of $\Sigma^*\times\Sigma^*$.
Unfortunately, the asynchronous automata introduced in \cite{RabinScott} do not behave well under function composition.
Therefore, we will use the following definition; it is a restricted form of the \emph{finite state transducers} used in computational linguistics.
Recall that we can view an automaton as a directed graph, with arrows between states labelled by input characters.
If we also add output words to these labels, we get a \emph{DOLT}: 

\begin{definition}
\vbox{
A \emph{DOLT} or \emph{DOmain-Led Transducer} $D$ consists of the data $(S,\Sigma,\mu,\rho,s_0)$.
Here:
\begin{itemize}
\item $S$ is a finite \emph{state set};
\item $\Sigma$ is a finite \emph{alphabet};
\item $s_0\in S$ is the \emph{initial state};
\item $\mu:S\times(\Sigma\amalg\{\$\})\partialto S$, a partial function, is the \emph{transition function}; and
\item $\rho:S\times(\Sigma\amalg\{\$\})\partialto(\Sigma\amalg\{\$\})^*$, with the same domain as $\mu$, is the \emph{output function}.
\end{itemize}}
We require that $\rho$ maps $S\times\Sigma$ into $\Sigma^*$ and $S\times\{\$\}$ into $\Sigma^*\$$, the set of words where $\$$ occurs exactly once, as the last letter.
As with finite state automata, the transition function induces $\mu:S\times(\Sigma\cup\{\$\})^*\partialto S$.
Likewise, we extend $\rho$ to $\rho:S\times(\Sigma\amalg\{\$\})^*\to(\Sigma\amalg\{\$\})^*$ by the recursion $\rho(s,uv)=\rho(s,u)\rho(\mu(s,u),v)$ where $\rho(s,\emptystring)=\emptystring$.  That is, $\rho(s,w)$ is the concatenation of the outputs as we follow the path labelled $w$ from $s$.
Set $L(D)=\{(x,y)\in\Sigma^*\times\Sigma^*\,|\,\rho(s_0,x\$)=y\$\}$.
\end{definition}

Note that $L(D)$ is the graph of a partial function.
Composition of relations will be defined by: $R\circ S=\{(u,v)\,|\,\exists w.(u,w)\in R\land(w,v)\in S\}$.
\textbf{Caution}: note the composition order commonly used for transducers.

\begin{lemma}\label{DOLT composition}
Given DOLTs $D_1,D_2$ with alphabet $\Sigma$ there is an effective procedure producing a DOLT $D_1\circ D_2$ satisfying $L(D_1\circ D_2)=L(D_1)\circ L(D_2)$.
\end{lemma}
\begin{proof}
We are given $D_i$ with state set $S_i$, transition function $\mu_i$, initial state $s_{0,i}$, and output function $\rho_i$.
For $D_1\circ D_2$, take state set $S=S_1\times S_2$. 
The transition function is $\mu((s_1,s_2),c)=(\mu_1(s_1,c),\mu_2(s_2,\rho_1(s_1,c)))$ with $c\in\Sigma\amalg\{\$\}$.
The initial state is $s_0=(s_{0,1},s_{0,2})$ and the output function is $\rho((s_1,s_2),c)=\rho_2(s_2,\rho_1(s_1,c))$.
\end{proof}

\begin{example}
Consider the DOLT $D$ with states $\{s_0,s_1,s_2\}$, where $(\mu,\rho)$ maps:
$$(s_0,a)\mapsto(s_0,\emptystring),\quad
 (s_0,b)\mapsto(s_1,a^2),\quad
 (s_1,a)\mapsto(s_1,b^2),\quad
 (s_1,b)\mapsto(s_1,a^2),\quad
 (s_1,\$)\mapsto(s_2,a\$).$$
Note that the fixed point set $\fix(L(D))=\{a^{2^{2n}}b^{2^{2n-1}}a^{2^{2n-2}}\cdots b^2a\,|\,n\geq1\}$ is not regular.
\end{example}

The remainder of this section is devoted to finding a suitable criterion on $D$ making $\fix(L(D))$ regular.
If $u,v,w\in\Sigma^*$ then $u$ is a \emph{prefix}, $w$ a \emph{suffix}, and $v$ a \emph{substring} of the word $uvw$.
Order prefixes by inclusion: $u\subseteq uv$.
Any word $w\in\Sigma^*$ decomposes uniquely into \emph{blocks}, maximal substrings of the form $c^n$ with $c\in\Sigma$.
Let $w[[i]]$ denote the prefix consisting of the first $i$ blocks.
(If $i$ is larger than the number of blocks, $w[[i]]=w$. If $i\leq0$, $w[[i]]=\emptystring$.)
If $w\in\dom(L(D))$, we say the \emph{image of $w[[i]]$} is $\rho(w[[i]]):=\rho(s_0,w[[i]])$.
A DOLT has \emph{Property $\mathcal{P}(N)$} if for all $(u,v)\in L(D)$ and all $i\in\mathbb{Z}$,
$v[[i-N]]\subseteq\rho(u[[i]])\subseteq v[[i+N]]$.
It has \emph{Property $\mathcal{P}(N,B)$} if, moreover, each block of each $w\in\dom(L(D))$ has length at most $B$.

\begin{example}
In the previous example, the word $a^{2^{2n}}b^{2^{2n-1}}a^{2^{2n-2}}\cdots b^2a$ decomposes into blocks
$a^{2^{2n}}$, $b^{2^{2n-1}}$, $a^{2^{2n-2}}$, $\ldots$ $b^2$, $a$.
The image of $a^{2^{2n}}$ is $\emptystring$, the image of $b^{2^{2n-1}}$ is $a^{2^{2n}}$, the image of $b^2$ is $a^4$, etc.
This DOLT enjoys Property $\mathcal{P}(1)$, but not $\mathcal{P}(1,B)$ for any $B$.
\end{example}

\begin{lemma}\label{fix(D) regular}
If DOLT $D$ enjoys Property $\mathcal{P}(N,B)$ then $\fix(L(D))$ is regular.
An automaton accepting $\fix(L(D))$ can be built effectively from $D$, $N$, and $B$.
\end{lemma}
\begin{proof}
We build an automaton $M$ which accepts $\fix(L(D))$.
The state $\mu_M(s_{0,M},u)=(\sigma,\widehat{u},\widehat{v})$ encodes:
\vbox{\begin{itemize}
\item the state $\sigma=\mu_D(s_{0,D},u)$,
\item a certain suffix $\widehat{u}$ of $u$, and
\item a certain suffix $\widehat{v}$ of $v:=\rho(s_{0,D},u)$.
\end{itemize}}
The intuition here is that these suffixes encode the portions of the input and output we intend to compare in the future, after comparing as much as possible in the present.
Thus at least one of $\widehat{u},\widehat{v}$ is always the empty string $\emptystring$.
The initial state $\mu_{0,M}=(\mu_{0,D},\emptystring,\emptystring)$.
The accept states of $M$ are those for which $\widehat{u}\$=\rho(\sigma,\$)$ and $\widehat{v}=\emptystring$ simultaneously.

Next we construct the data for $uc$, with $c\in\Sigma$, from the data for $u$.
Let $\widehat{u}'$ and $\widehat{v}'$ be obtained from $\widehat{u}c$ and $\widehat{v}\rho(\sigma,c)$, respectively, after deleting their longest common prefix.
If $\widehat{u}'$ and $\widehat{v}'$ are both nonempty, or if $\widehat{v}'$ contains a block of length $>B$,
then $uc$ cannot be a prefix of a fixed point.
Otherwise, we give $M$ a transition arrow from $(\sigma,\widehat{u},\widehat{v})$ to $(\mu_D(\sigma,c),\widehat{u}',\widehat{v}')$
labelled $c$.

It remains to show that $M$ requires boundedly many states.
Suppose $u$ is a prefix of fixed point $w$.
Then $u$ contains at most $N$ more blocks than $v:=\rho(u)$.
Therefore $\widehat{u}$ has length at most $(N+1)B$.
Similarly, $\widehat{v}$ has length at most $(N+1)B$.  Therefore we can take finite state set:
$$S_M=\left\{(\sigma,\widehat{u},\widehat{v})\in S_D\times\Sigma^*\times\Sigma^*\,:\,
     |\widehat{u}|,|\widehat{v}|\leq (N+1)B,\,
		 |\widehat{u}||\widehat{v}|=0\right\}.$$
Thus $M$ is a finite state automaton accepting $\fix(L(D))$.
\end{proof}

\begin{lemma}\label{property P composition}
If DOLT $D_1$ has Property $\mathcal{P}(N_1)$ and $D_2$ has $\mathcal{P}(N_2)$ then $D_1\circ D_2$ has $\mathcal{P}(N_1+N_2)$.
\end{lemma}
\begin{proof}
Suppose $(u,v)\in L(D_1)$ and $(v,w)\in L(D_2)$.
Then $(\rho_1\circ\rho_2)(u[[i]])\subseteq\rho_2(v[[i+N_1]])\subseteq w[[i+N_1+N_2]]$.
Likewise, $(\rho_1\circ\rho_2)(u[[i]])\supseteq w[[i-(N_1+N_2)]]$.
\end{proof}

\begin{example}\label{Example:free group multiplication}
Consider a finite rank free group $F(X)$.
Let $\Sigma=X\amalg\ol{X}$ with $\ol{X}$ a set of formal inverses. 
Each element $w\in F(X)$ is represented by a unique reduced word $\redu{w}\in\Sigma^*$ of length $|w|$.
Fix $u,v\in F(X)$.
Let $f_{u,v}$ be the (partial) function $\Sigma^*\partialto\Sigma^*$ with graph $\{(\redu{w},\redu{uwv}):w\in F(X)\}$.
There is a DOLT with Property $\mathcal{P}(|u|+|v|)$ accepting $f_{u,v}$.
The construction from $u,v$ is effective.
\end{example}
\begin{proof}
Since $f_{u,v}=f_{u,\emptystring}\circ f_{\emptystring,v}$, it suffices by Lemma \ref{property P composition}
       to assume $u=\emptystring$ or $v=\emptystring$.
We do the case $u=\emptystring$.

Let $D=\Sigma^{\leq|v|}$.
The initial state is $s_{0}:=\emptystring$.
Suppose we are in state $\alpha$ and we encounter input letter $c\in\Sigma$.
Assume $\alpha c$ is reduced (otherwise, there will be no transition arrow).
If $|\alpha c|<|v|$, we transition to state $\alpha c$ and output $\emptystring$.
Otherwise, write $\alpha c=\beta\gamma$ with $\beta\in\Sigma$ and $\gamma\in\Sigma^*$; 
           transition to state $\gamma$ and output $\beta$.

If we are in state $\alpha$ and encounter input $\$$, output $\redu{\alpha v}\$$.

By construction, $L(D)=f_{\emptystring,v}$ and $D$ enjoys Property $\mathcal{P}(|v|)$.
The case $v=\emptystring$ is very similar.
\end{proof}

\section{The conjugacy problem for Higman's group}\label{Section:Equation}
In this section, we show the conjugacy problem in $H$ is decidable.
As indicated in section \ref{Section:Warmup},
  the starting point is an analysis of the equation $\gamma_1a_1=a_1'\gamma_2$ where $a_1,a_1'\in A\setminus C$ are given.
Let
	$$D(a_1,a_1')=\{(\gamma_1,\gamma_2)\in\{a^\pm,c^\pm\}^* \,|\, \gamma_1a_1=_Aa_1'\gamma_2\textrm{ and $\gamma_i$ is freely reduced}\}.$$
Recall from section \ref{Section:Automata} our composition convention 
                                           and the partial function $f_{u,v}:\{a^\pm,c^\pm\}^*\partialto\{a^\pm,c^\pm\}^*$.

\begin{lemma}\label{Lemma:reduction to tight}
If $a_1=_Au\widehat{a_1}v$ and $a_1'=_AU\widehat{a_1}'V$ with $u,v,U,V\in C$ then:
$$D(a_1,a_1')=f_{U^{-1},u}\circ D(\widehat{a_1},\widehat{a_1}')\circ f_{V^{-1},v}.$$
\end{lemma}
\begin{proof}
$\gamma_1a_1=_Aa_1'\gamma_2
\iff U^{-1}\gamma_1u\widehat{a_1}=_A\widehat{a_1}'V\gamma_2v^{-1}
\iff(\redu{U^{-1}\gamma_1u},\redu{V\gamma_2v^{-1}})\in D(\widehat{a_1},\widehat{a_1}').$
\end{proof}

Lemma \ref{Lemma:reduction to tight} lets us restrict attention to \emph{tight} $a_1,a_2$ (defined below) when solving $\gamma_1a_1=a_1'\gamma_2$.

\begin{definition}
Consider a Britton--reduced word $w=(r_1,s_1)a^{n_1}(r_2,s_2)a^{n_2}\cdots a^{n_{k-1}}(r_k,s_k)$ 
      in $\{a^+,a^-\}\amalg\mathbb{Z}[1/2]\rtimes\mathbb{Z}$ representing an element of $A\setminus C$, with $k\geq1$.
We say $w$ can be \emph{tightened} if $w=_Auw'v$ for some $u,v\in C$ and $w'=(r_1',s_1')a^{n_1'}\cdots a^{n_{j-1}'}(r_j',s_j')$ and: $j<k$ OR
$w'$ (but not $w$) is an odd power of $b$.
Otherwise, $w$ is \emph{tight}.
\end{definition}

\begin{lemma}\label{Lemma:factorize a1}
Given a Britton--reduced word $w$ in $\{a^+,a^-\}\amalg\mathbb{Z}[1/2]\rtimes\mathbb{Z}$ representing an element of $A\setminus C$,
we can effectively produce a factorization $w=_Auw'v$, with $w'$ tight and $u,v\in C$.
\end{lemma}

\begin{proof}
Tightening $w=(r_1,s_1)a^{n_1}(r_2,s_2)a^{n_2}\cdots a^{n_{k-1}}(r_k,s_k)$, $k\geq2$, on the left side is possible if and only if:
${r_1\in\mathbb{Z}}$ and $a^{-\sgn(n_1)}c^{-r_1}w$ is \emph{not} Britton--reduced.
Tightening on the right is similarly effective.
Now, ${w=(r_1,s_1)}$, $r_1\neq0$ can be tightened if and only if $r:=r_1/(1-2^{-s_1})\in\mathbb{Z}$.
In this case, $c^{-r}wc^r$ is of the form $(0,s)$.
Finally, to tighten $w=(0,s_1)$, write $s_1=2^nm$ with $m$ odd; then $w=a^{-n}b^ma^n$.
\end{proof}

\begin{lemma}\label{one candidate}
Suppose $a_1=(r_1,s_1)a^{n_1}(r_2,s_2)a^{n_2}\cdots a^{n_{k-1}}(r_k,s_k)$ and $a_1'=(r'_1,s'_1)\cdots$ are tight and $k\geq2$.
Then $$\dom D(a_1,a_1')\subseteq\{c^{r'_1-r_1}\}.$$
In particular, we can effectively build a Property $\mathcal{P}(1,|r_1'-r_1|)$ DOLT accepting $D(a_1,a_1')$.
\end{lemma}
\begin{proof}
Consider the Cayley 2--complex for the presentation $$A=\left\langle a,b,c \, \left| \,  b^a = b^2, c^b = c^2 \right. \right\rangle.$$
In any reduced van Kampen diagram, cells with interior edges labelled $a$ join up along these edges into \emph{$a$-corridors}.
(There can be no $a$--annuli since $\langle b\rangle$ is torsion--free.)

Suppose $(\gamma_1,\gamma_2)\in D(a_1,a_1')$.
Consider a reduced van Kampen diagram establishing $\gamma_1a_1=a_1'\gamma_2$.
There must be an $a$-corridor $\mathcal{C}$ originating on side $a_1$ because $k\geq2$; consider an innermost such corridor.
The other end of $\mathcal{C}$ cannot be on side $a_1$ (because $a_1$ is Britton--reduced)
                                     or on sides $\gamma_1,\gamma_2$ (because $a_1$ is tight).
Therefore, $\mathcal{C}$ joins $a_1$ to $a_2$.

So there can be no $a$-corridor between $\gamma_1$ and $\gamma_2$.
Now, $\gamma_1$ is Britton--reduced since $\langle c\rangle\cap\langle b\rangle=\{1\}$.  So $\gamma_1=c^j$ for some $j\in\mathbb{Z}$.
Considering an $a$-corridor originating at $a^{n_1}$ on $a_1$, we see
$$c^j(r_1,s_1)\in(r_1',s_1')\langle b\rangle.$$
Hence $j=r_1'-r_1$, so $\gamma_1=c^{r_1'-r_1}$.
(If $r_1'-r_1\notin\mathbb{Z}$, we can immediately conclude $D(a_1,a_1')=\emptyset$.)

Note that checking whether $c^{r_1'-r_1}\in\dom D(a_1,a_1')$ and finding its output comes down to 
checking whether $a_1'^{-1}c^{r_1'-r}a_1\in C$ and, if so, writing it in $\{a^\pm,c^\pm\}^*$.
So building the desired DOLT recognizing $D(a_1,a_1')$ is effective by Lemma \ref{Britton effective}.
\end{proof}

\vbox{
\begin{lemma}\label{affine case}
Suppose $a_1=(r,s)$ is tight and not in $\langle b\rangle$.
If $a_1'$ is tight and $D(a_1,a_1')\neq\emptyset$ then $a_1'=(r',s)\notin\langle b\rangle$.
Write $r-r'=\frac{p}{2^q}$ with $p$ odd (take $q=-\infty$ if $r-r'=0$).  We have cases:
\begin{itemize}
\item if $q>\max\{s,0\}$ then $D(a_1,a_1')=\emptyset$;
\item if $s\geq\max\{q,0\}$ then $D(a_1,a_1')=\{(c^n,c^{f(n)}):n\in\mathbb{Z}\}$ with $f(n)=2^s(r-r'+n)$;
\item if $s,q\leq0$ then $D(a_1,a_1')=\{(c^{g(m)},c^m):m\in\mathbb{Z}\}$ with $g(m)=2^{-s}m+r'-r$.
\end{itemize}
In particular, $D(a_1,a_1')$ is accepted by a Property $\mathcal{P}(1)$ DOLT effectively constructed from $r,r',s$.
\end{lemma}}

\begin{proof}
Consider a solution $\gamma_1(r,s)=a_1'\gamma_2$.  By tightness, there is no $a$-corridor in a reduced van Kampen diagram.
Therefore $\gamma_1=c^n$ and $\gamma_2=c^m$ and $a_1=(r',s')$ for some $n,m,s'\in\mathbb{Z}$ and $r'\in\mathbb{Z}[1/2]$.
The equation
$$\gamma_1(r,s)=a_1'\gamma_2\iff(n,0)(r,s)=(r',s')(m,0)\iff(n+r,s)=(r'+2^{-s'}m,s')$$
is equivalent to $s'=s$ and $m=2^s(r-r'+n)=2^s\left(\frac{p}{2^q}+n\right)$.
The cases follow immediately.
\end{proof}

\begin{definition}
For $w\in\{a^\pm,c^\pm\}^*$, the \emph{$a$--height} of a letter is the $a$--exponent sum of the prefix ending on that letter.
An \emph{$a$--Dyck word} $w$ is a freely reduced word in $\{a^\pm,c^\pm\}^*$ with $a$--exponent sum zero,
so that each letter has non--negative $a$--height.
Define \emph{$c$--height} and \emph{$c$--Dyck} analogously.
\end{definition}

\begin{example} In the $a$--Dyck word $c^{-2}a^2c^3a^{-1}c^4ac^2a^{-2}cac^{-3}a^{-1}c^4$,
the blocks $c^{-2}$, $c^3$, $c^4$, $c^2$, $c$, $c^{-3}$, $c^4$ occur at $a$--heights 0, 2, 1, 2, 0, 1, 0 respectively.
\end{example}

\begin{lemma}\label{Dyck case}
Suppose $s,s'$ are odd.
Then $D(b^s,b^{s'})=\emptyset$ unless $s=s'$.
Set $\mathcal{D}_{s}=D(b^{s},b^{s})$.

For $s>0$, the domain of $\mathcal{D}_{s}$ is the set of $a$--Dyck words.
The output of such a word is the new $a$--Dyck word obtained by replacing each $c^\pm$ at height $h$ by $c^{\pm2^{s\cdot 2^h}}$.
For $s<0$, we have $\mathcal{D}_{s}=\mathcal{D}_{-s}^{-1}$.

For any fixed $H,s$, we can effectively build a DOLT with Property $\mathcal{P}(0)$ which accepts 
$$\mathcal{D}_{s,H}:=\{(u,v)\in\mathcal{D}_s\,|\,\textrm{ $a$-height of each letter $\leq H$}\}.$$
\end{lemma}
\begin{example} $(c^5acac^{-1}a^{-2}c,c^{40}ac^{64}ac^{-4096}a^{-2}c^8)\in\mathcal{D}_3$;
 $\,\,\,(c^{40}ac^{64}ac^{-4096}a^{-2}c^8,c^5acac^{-1}a^{-2}c)\in\mathcal{D}_{-3}$.
\end{example}

\begin{proof}[Proof of Lemma \ref{Dyck case}]
Consider a reduced van Kampen diagram for $\gamma_1b^{s}=b^{s'}\gamma_2$.
Each $a$--corridor runs from $\gamma_1$ to $\gamma_2$.
Order the corridors by the orientation of $\gamma_1$, from $b^{s'}$ (back) to $b^s$ (front).
Each corridor contributes either $a^{+1}$ or $a^{-1}$ to $\gamma_1$, and the same letter to $\gamma_2$.
The other sides of the corridor are labelled, in some order, $b^i$ and $b^{2i}$ for some $i\in\mathbb{Z}$.
(The front side is $b^{2i}$ if and only if the letter contributed to $\gamma_1,\gamma_2$ is $a^{+1}$.)

Besides the $a$--corridors, each 2--cell has label $c^bc^{-2}$, with $b$--exponent sum zero.
So the back--most $a$--corridor must have back side labelled $b^{s'}$; the front--most corridor must have front side labelled $b^s$;
 and consecutive corridors must share a label.
In particular, each label $b^i$ considered shares the same odd part.  So $s'=s$.

By similar reasoning, the $a$--exponent sums on $\gamma_1,\gamma_2$ are zero
   (there are as many doubling corridors as halving corridors since $s'=s$).
Further, for each $a$--corridor, the letter $a^\pm$ contributed to each of $\gamma_1,\gamma_2$ has height$\geq0$.
Otherwise, its front edge label would not be an integral power of $b$.
Therefore, $\gamma_1,\gamma_2$ are $a$--Dyck words.

Conversely, for $s>0$, an easy calculation shows that each $a$--Dyck $\gamma_1$ produces the value $\gamma_2:=b^{-s}\gamma_1b^s$
claimed in the Lemma's statement.  The fact $\mathcal{D}_{-s}=\mathcal{D}_s^{-1}$ follows by reflecting the van Kampen diagrams.

It remains to construct the DOLT for arbitrary $H,s$.
In the case $s>0$, the machine need only keep track of the $a$--exponent sum.
(Additional states are needed to ensure $\gamma_1$ is freely reduced, but we ignore these for simplicity.)
Take state set $\Sigma=\{0,\ldots,H\}$.
We have $(\mu,\rho)(h,a^\pm)=(h\pm1,a^\pm)$, if $h\pm1\in\Sigma$, and $(\mu,\rho)(h,c^\pm)=(h,c^{\pm2^{s\cdot2^h}})$.
Finally, $\rho(0,\$)=\$$.  Thus we only accept words with $a$--exponent sum zero.

We merely indicate the changes needed for the case $s<0$.
When the DOLT is at $a$--height $h$, it should only output $c^\pm$ once at the end of every $2^{s\cdot2^h}$ instances of input $c^\pm$.
So the DOLT needs states tracking how long it has been since such an output.
If fed $a^\pm$ or $\$$ in the middle of such a cycle, the DOLT should fail.

Finally, Property $\mathcal{P}(0)$ follows from the fact that $\mathcal{D}_s$ (or $\mathcal{D}_s^{-1}$) is computed by
   $a^\pm\mapsto a^\pm$, $c^\pm\mapsto c^{\pm2^{s\cdot2^h}}$.
\end{proof}

\vbox{In summary, given tight $a_1,a_1'\in A\setminus C$, we can determine whether $\gamma_1a_1=a_1'\gamma_2$ solutions:
\begin{itemize}
\item are unique (Lemma \ref{one candidate}),
\item involve $c$--blocks (Lemma \ref{affine case}),
\item involve $a$--Dyck words (Lemma \ref{Dyck case}), or
\item do not exist.
\end{itemize}}

We are also interested in solving $\gamma_2b_2=b_2'\gamma_3$ given $b_2,b_2'\in B\setminus C$.
Note that $A\cong B$ via an isomorphism sending $a\mapsto c, b\mapsto d, c\mapsto a$ (and thus preserving $C$).
So Lemmas \ref{Lemma:reduction to tight}, \ref{Lemma:factorize a1}, \ref{one candidate}, \ref{affine case}, \ref{Dyck case} still apply;
we just reverse the roles of $a$ and $c$.
In particular, solutions $(\gamma_2,\gamma_3)\in C^2$ to $\gamma_2b_2=b_2'\gamma_3$ can involve $a$--blocks and $c$--Dyck words.
Write $\mathcal{D}'_s$ for the analogue to $\mathcal{D}_s$ solving $\gamma_2d^s=_Bd^s\gamma_3$.
Define $\mathcal{D}'_{s,H}$ similarly.

\begin{example} $(a^{40}ca^{64}ca^{-4096}c^{-2}a^8,a^5caca^{-1}c^{-2}a)\in\mathcal{D}'_{-3}$. \end{example}

\begin{lemma}\label{C not conjugate to F}
Suppose $x\in C=\langle a,c\rangle$ and $y\in F=\langle b,d\rangle$ are conjugate in $H$.  Then $x=1$.
\end{lemma}
\begin{proof}
We may assume $y$ is cyclically reduced as an element of $F$.
By Lemma \ref{conjugation for amalgams}, $x$ must have the same length as $y$ as a word in $A\amalg B$.
Therefore $y$ is a power of $b$ or of $d$.  Without loss of generality, assume $y=b^n$.
Considering the Britton--reduced form of $z$ making $x=y^z$, we conclude $y$ is conjugate in $A$ to an element of $C$.
Therefore we may assume $x\sim_Ay=b^n$.

Now, just repeat the argument:
We may assume $x$ is cyclically reduced as an element of $C$.
Consider the decomposition $A=\left\langle a,b \, \left| \,  b^a = b^2 \right. \right\rangle
       \ast_{\langle b\rangle}\left\langle b,c \, \left| \,  c^b = c^2 \right.\right\rangle$.
This time, Lemma \ref{conjugation for amalgams} implies $x$ lies in one of the two factors.
Without loss of generality, $x=c^m$ and $(m,0)\sim_K(0,n)$.  Therefore $n=0$ so $x=y=1$.
\end{proof}

We record the following easy fact:
\begin{lemma}\label{easy fact}(\cite[Theorem 4.6]{MKS})
Suppose $x,y\in H$ are cyclically Britton--reduced, and that $x\in A$ is not conjugate in $H$ to any element of $C$.
Then $x\sim_Hy$ if and only if $x\sim_Ay$.
\end{lemma}

\begin{proof}[Proof of Theorem \ref{conjugacy in H}]
We are given words $x,y\in\{a^\pm,b^\pm,c^\pm,d^\pm\}^*$ and must determine whether $x\sim_Hy$.

By Lemma \ref{Britton effective}, we can replace $x,y$ by cyclically Britton--reduced words in $A\amalg B$.
We now reduce to the case that $x$ has length at least 2, in order to apply Lemma \ref{conjugation for amalgams}.
Let $\psi:H\to H$ denote the automorphism mapping $a\mapsto b\mapsto c\mapsto d\mapsto a$.
Let $x',y'$ be cyclic Britton--reductions of $\psi(x),\psi(y)$.
If any of $x,y,x',y'$ have length at least 2, the reduction is complete (since $x\sim_Hy\iff x'\sim_Hy'$).
So suppose they all have length 1.

For simplicity of notation, suppose that (say) $x\in A$ and $x'\in B$.
Now, $x$ and $x'$ cannot both be conjugate in $H$ to elements of $C$:
if they were, then $x$ would be conjugate to an element of $\psi^{-1}(C)=F$, contradicting Lemma \ref{C not conjugate to F}.
Therefore $x\sim_Hy$ if and only if $x\sim_Ay$ or $x'\sim_By'$ (Lemma \ref{easy fact}).
Since the conjugacy problem in $A$ is decidable (Proposition \ref{conjugacy in A}), we are done in this case.

We have reduced to the case $x=a_1b_2\cdots a_{2n-1}b_{2n}$, for some $n\geq1$, with $a_i\in A$ and $b_j\in B$.
If $y$ has different length from $x$, we immediately conclude that $x\nsim_Hy$.
Thus $y=a_1'b_2'\cdots a_{2n-1}'b_{2n}'$, with $a_i'\in A$, $b_j'\in B$.
By Lemma \ref{conjugation for amalgams}, we need only check if $zxz^{-1}=y'$ for some $z\in C$ and some cyclic permutation $y'$ of $y$.
These permutations can be checked one at a time.

By Lemmas \ref{Lemma:reduction to tight}, \ref{Lemma:factorize a1}, \ref{one candidate}, \ref{affine case}, and \ref{Dyck case},
we are reduced to determining whether
$$\Phi:=f_{u_0,v_0}\circ D_1\circ f_{u_1,v_1}\circ D_2\circ f_{u_2,v_2}\circ\cdots\circ D_{2n}\circ f_{u_{2n},v_{2n}}$$
has a fixed point in $C$.
Here each $u_i,v_i$ is an effectively constructed element of $C$ and each $D_i$ is either an effectively constructed DOLT
or $\mathcal{D}_s$ or $\mathcal{D}'_s$ for some algorithmically determined $s$.
We are free to cyclically permute $\Phi$; this does not affect existence of fixed points.
In particular, we may absorb $f_{u_0,v_0}$ into $f_{u_{2n},v_{2n}}$.

At this point, we do not know that $\Phi$ is accepted by a DOLT.
Nonetheless, it is the graph of a computable function by Lemma \ref{Dyck case}.
If we can find a finite set $S$ containing the range of some $D_i$ or some $D_i\circ f_{u_i,v_i}\circ D_{i+1}$, 
then we can compute a finite set $T$ containing the range of $\Phi$.
This would solve the conjugacy problem, as each $z\in T$ can be tested to see if $zxz^{-1}=_Hy'$.

In particular, we may assume no $D_i$ is a singleton (as in Lemma \ref{one candidate}).
Suppose next that $D_i$ accepts $a$--blocks and $D_{i+1}$ accepts $c$--blocks (as in Lemma \ref{affine case}).
Since outputs of $D_i\circ f_{u_i,v_i}$ contain at most $|u_i|+|v_i|$ instances of $c^\pm$, Lemma \ref{affine case}
bounds the range of $D_i\circ f_{u_i,v_i}\circ D_{i+1}$. This rules out consecutive block--case DOLTs.

Next, suppose $D_i$ accepts $a$--blocks and $D_{i+1}=\mathcal{D}_s$ works with $a$--Dyck words.
The outputs of $D_i$ with length $>|u_i|+|v_i|$ yield outputs of $D_i\circ f_{u_i,v_i}$ with non-zero $a$--exponent sum.
These cannot be valid inputs to $D_{i+1}$.
We have again bounded the range of $\Phi$, so block--type DOLTs are ruled out.

At this point, we have each $D_i$ of the form $\mathcal{D}_s$ or $\mathcal{D}'_s$.
We wish to replace each $\mathcal{D}_s$ in the definition of $\Phi$ by $\mathcal{D}_{s,H}$ (and $\mathcal{D}'_{s'}$ by $\mathcal{D}'_{s',H}$)
 to get $\Phi_H$, accepted by some DOLT $D$, with $\fix(\Phi_H)=\fix(\Phi)$.
Take $N$ larger than all $|u_i|+|v_i|$.
Now choose $H$ satisfying $2^{2^H}\geq N+H$ (for example, $H=N$).
In the following, we use \emph{height} to mean $a$--height for $D_i=\mathcal{D}_s$ and $c$--height for $D_i=\mathcal{D}'_{s'}$.
Suppose, for a contradiction, that $(w,w)\in\Phi\setminus\Phi_H$.
Then some $D_i$ achieves height $h>H$ on its corresponding input $u$, with output $v$.
Take $i,h$ to make $h$ maximal (with $w$ fixed).
For simplicity, assume $D_i=\mathcal{D}_s$.
Then one of $u,v$ contains a $c^\pm$--block of size $\geq 2^{2^h}$, by Lemma \ref{Dyck case}.
The corresponding input/output to $D_{i\pm1}$, indices taken cyclically, contains a $c^\pm$--block of size $\geq 2^{2^h}-N>h$.
This contradicts the maximality of $h$, so $\fix(\Phi)=\fix(\Phi_H)$.

Finally, we show $D$ enjoys property $\mathcal{P}(N,B)$ for effectively constructible $N,B$.
Take $N=\sum_i{|u_i|+|v_i|}$.
Each $\mathcal{D}_s,\mathcal{D}'_{s'}$ enjoys Property $\mathcal{P}(0)$ by Lemma \ref{Dyck case}.
So we get Property $\mathcal{P}(N)$ for $D$ from Lemma \ref{property P composition} \& Example \ref{Example:free group multiplication}.
Suppose $D_1=\mathcal{D}_{s,H}$ and $D_2=\mathcal{D}'_{s',H}$.
The length of an $a^\pm$ input block for $D_1$ is bounded by $H$.
Similarly, the length of a $c^\pm$ output block accepted by $f_{u_1,v_1}\circ D_2$ is at most $H+|u_1|+|v_1|$.
The corresponding input $c^\pm$ block for $D_1$ is at most $B:=2^{|s|\cdot 2^{H+|u_1|+|v_1|}}\geq H$, by Lemma \ref{Dyck case}.
So $D$ enjoys property $\mathcal{P}(N,B)$.

Therefore Lemma \ref{fix(D) regular} shows $\fix(\Phi)$ is accepted by an effectively constructed automaton $M$.
Lemma \ref{emptiness problem} then tells us whether $\fix(\Phi)$ is empty.  This solves the conjugacy problem.
\end{proof}

\section{Generic Case Complexity of the Conjugacy Problem} \label{Section:Runtime}

In this section, we establish a time $\bigo(n^7)$ solution
to the conjugacy problem for $H$ in a strongly generic setting (Proposition \ref{H conjugacy runtime} and
    Theorem \ref{strongly generic conjugacy in H}).
To decide whether $x\sim_Hy$, we examine cyclic Britton--reductions $\widehat{x}$ and $\widehat{y}$.
As we will see (Lemmas \ref{nice with backtracking} and \ref{nice without backtracking}), $\widehat{x}$ or $\widehat{y}$ typically contains a \emph{nice factor} (definition below) -- an entry to which Lemma \ref{one candidate} applies after tightening.
This leads to a unique candidate conjugator, which can be checked using a solution to the word problem.
The word problem in $H$ is solved in time $\bigo(n^6)$ by \cite{DLU2}, using \emph{power circuit} technology introduced in \cite{MUW}.

We recall the definition of power circuit.
Let $\Gamma$ be a finite set and $\delta:\Gamma\times\Gamma\to\{-1,0,1\}$ a map.
The \emph{support} $\{(P,Q)\in\Gamma\times\Gamma:\delta(P,Q)\neq0\}$ is required to be a directed acyclic graph.
(In particular, $\delta(P,P)=0$.)
A \emph{marking} is a mapping $M:\Gamma\to\{-1,0,1\}$.
It has \emph{support} $\sigma(M)=\{P\in\Gamma:M(P)\neq0\}$.
For each $P\in\Gamma$, there is an associated marking $\Lambda_P(Q):=\delta(P,Q)$.
\emph{Evaluation} of markings $\evale(M)$ and of nodes $\evale(P)$ is defined by simultaneous recursion:
$$\evale(\emptyset)=0,\qquad\evale(P)=2^{\evale(\Lambda_P)},\qquad\evale(M)=\sum_P{M(P)\evale(P)}.$$
If each node (equivalently, marking) evaluates to an integer, $(\Gamma,\delta)$ is called a \emph{power circuit} of size $|\Gamma|$.
\vbox{A \emph{reduced power circuit} $\Pi$ consists of:
\begin{itemize}
\item a power circuit $(\Gamma,\delta)$ in which no two nodes evaluate to the same integer;
\item further data used only in a ``black box'' algorithm, \texttt{ExtendReduction}, within \cite{DLU2}.
Without this improvement of \cite{DLU2} over \cite{DLU1}, there would be additional logarithmic factors in Lemma \ref{s steps} below.
\end{itemize}}

An element $(r,s)\in\mathbb{Z}[1/2]\rtimes\mathbb{Z}$ can be represented (non--uniquely) by a triple of integers
$$[u,v,w]:=(0,v)(u,w)=(2^{-v}u,v+w)$$ where $u,v,w\in\mathbb{Z}$ and $v\geq0\geq w$.
A \emph{triple marking} $[U,V,W]$ consists of markings $U,V,W$ on a single reduced power circuit,
representing $[\evale(U),\evale(V),\evale(W)]\in\mathbb{Z}[1/2]\rtimes\mathbb{Z}$.
Equipping a triple marking with a \emph{type} among $G_{ab}, G_{bc}, G_{cd}, G_{da}$ gives it a semantic value in $H$ via
\begin{align*}
H=A\ast_CB
&=(G_{ab}\ast_{\langle b\rangle}G_{bc})\ast_{\langle a,c\rangle}(G_{cd}\ast_{\langle d\rangle} G_{da})\\
&=\left(\left\langle a,b\,\left|\, b^a=b^2\right.\right\rangle \ast_{\langle b\rangle}
        \left\langle b,c\,\left|\, c^b=c^2\right.\right\rangle\right)
	\ast_{\langle a,c\rangle}
	\left(\left\langle c,d\,\left|\, d^c=d^2\right.\right\rangle \ast_{\langle d\rangle}
				\left\langle d,a\,\left|\, a^d=a^2\right.\right\rangle\right).
\end{align*}

In \cite{DLU2}, an element of $H$ is represented by a \emph{main data structure} $\mathcal{T}=(\Pi,(T_j)_{j\in J})$,
where $\Pi$ is a reduced power circuit and $(T_j)_{j\in J}$ is a finite sequence of triple markings on $(\Gamma,\delta)$ of various types.
The sequence decomposes into \emph{intervals}:
intervals of type $A$ are maximal subsequences with types among $G_{ab}$ and $G_{bc}$;
intervals of type $B$ are maximal subsequences with types among $G_{cd}$ and $G_{da}$.
Thus $\mathcal{T}$ is a \emph{power circuit representation} of both a word in $(A\amalg B)^*$ and an element of $H$.

\begin{definition}
Each triple marking $[U,V,W]$ has \emph{weight} $\omega([U,V,W])=|\sigma(U)|+|\sigma(V)|+|\sigma(W)|$.
The \emph{weight} of a main data structure $\mathcal{T}$ is $\omega(\mathcal{T})=\sum_{j\in J}{\omega(T_j)}$.
The \emph{size} of $\mathcal{T}$ is $\|\mathcal{T}\|=|\Gamma|$.
\end{definition}

\vbox{To manipulate power circuit representations, \cite{DLU2} defines various \emph{basic operations}.
The ones for $A$ are listed below (those for $B$ are analogous):

\begin{itemize}
\item Multiplication: 
        $$[u,v,w]_{ab}\cdot[u',v',w']_{ab}=[2^{v'}u+2^{-w}u',v+v',w+w']_{ab};$$
				$$[u,v,w]_{ab}\cdot[u',v',w']_{bc}=[2^{v'}u+2^{-w}u',v+v',w+w']_{bc};$$
\item Swapping $G_{ab}$ to $G_{bc}$:
        $$[u,0,0]_{ab}=\begin{cases}[0,u,0]_{bc}\textrm{ if $u\geq0$}\\ [0,0,u]_{bc}\textrm{ if $u<0$}\end{cases};$$
\item Swapping $G_{bc}$ to $G_{ab}$: 
        $$[0,v,w]_{bc}=[v+w,0,0]_{ab};$$
\item Splitting as $\langle a\rangle\langle b\rangle$ and $\langle c\rangle\langle b\rangle$:
           \begin{align*}[u,v,w]_{ab}&=[2^{-v}u,0,0]_{ab}\cdot[0,v,w]_{ab}\qquad\textrm{ if $2^{-v}u\in\mathbb{Z}$},\\
           [u,v,w]_{bc}&=[0,v,w]_{bc}\cdot[2^wu,0,0]_{bc}\qquad\textrm{ if $2^{v}u\in\mathbb{Z}$}.\end{align*}
\end{itemize}}

To apply a basic operation, replace the left side with the right side and forget the replaced markings.
For our current purposes, we add:
\begin{itemize}
\item Splitting as $\langle b\rangle\langle a\rangle$ and $\langle b\rangle\langle c\rangle$: 
           \begin{align*}[u,v,w]_{ab}&=[2^{-v}u,0,0]_{ab}\cdot[0,v,w]_{ab}\qquad\textrm{ if $2^{-v}u\in\mathbb{Z}$},\\
           [u,v,w]_{bc}&=[0,v,w]_{bc}\cdot[2^wu,0,0]_{bc}\qquad\textrm{ if $2^{v}u\in\mathbb{Z}$}.\end{align*}
\end{itemize}

Even with the new operations, the verbatim proof of \cite[Prop.\ 18]{DLU2} gives:

\begin{lemma}\label{s steps} 
Suppose $\mathcal{T}=(\Pi,(T_j)_{j\in J})$ is a main data structure with size $\leq m$, weight $\leq w$, and $|J|+w\leq m$.
Performing a sequence of $s$ basic operations (including tests like $2^vu\in\mathbb{Z}$) takes time $\bigo(s^2m^2)$.
The weight $\omega(\mathcal{T})$ does not increase, and the size $\|\mathcal{T}\|$ remains bounded by $\bigo(m+sw)$.\qed
\end{lemma}

\vbox{\begin{definition}
Consider $a_i=u\widehat{a_i}v\in A-C$ with $u,v\in C$ and $\widehat{a_i}$ tight.
Suppose $\widehat{a_i}\in(\{a^+,a^-\}\amalg\mathbb{Z}[1/2]\rtimes\mathbb{Z})^*$ contains at least one occurrence of $a^+$ or $a^-$.
(Thus $\widehat{a_i}$ is precisely the type of word to which Lemma \ref{one candidate} applies.)
In this case, we say $a_i$ is \emph{nice}.
We also call \emph{nice} the image of $a_i$ in $B-C$, under the isomorphism $A\cong B:a\mapsto c, b\mapsto d, c\mapsto a$.
In the Britton--reduced word $a_1b_2\cdots a_{2n-1}b_{2n}\in(A\amalg B)^*$, with $a_i\in A, b_j\in B$ and $n\geq1$,
we call any nice $a_i$ or $b_j$ a \emph{nice factor}.
\end{definition}}

\begin{prop}\label{H conjugacy runtime}
The following can be computed in time $\bigo(n^7)$:
Input words $x,y\in\Sigma^*=\{a^\pm,b^\pm,c^\pm,d^\pm\}^*$ of total length $n$.
Determine cyclically Britton--reduced forms $\widehat{x},\widehat{y}$ and decide if either has length $\geq2$ (as a word in $(A\amalg B)^*$) and a nice factor.
If so, determine whether $x\sim_Hy$ and find a power circuit representation of $z$ with $zxz^{-1}=y$ (if one exists).
\end{prop}

\begin{proof}
The word problem is solved in \cite{DLU2} as follows. 
Given a word $x$ of length $n$ in $\{a^\pm,b^\pm,c^\pm,d^\pm\}$,
a main data structure of size $m=\bigo(|x|)$ and weight $w=\bigo(|x|)$ is constructed to store it.
Then $\bigo(|x|^2)$ basic operations are performed until 
      the main data structure holds a Britton--reduced sequence of intervals representing $x$.
This takes time $\bigo(|x|^6)$.
Cyclic Britton--reduction can be achieved with the same bounds on time and number of operations,
 by Lemma \ref{cyclic reduction amalgams}.
			\footnote{Cf. \cite[Prop. 4]{DMW2}, where the analogous observation is made for Britton--reduction in $G_{1,2}$.}
Thus we have cyclically Britton--reduced $\widehat{x}$ and $\widehat{y}$ using $\bigo(n^2)$ basic operations.

We now test each interval of $\widehat{x}$ and $\widehat{y}$ to see if it represents a nice factor.
Without loss of generality, consider an interval $\mathcal{L}$ of type $A$ and sequence length $s\leq n$.
Recall that tightening and Britton--reduction in $A$ were defined viewing $A$ as an HNN--extension over $K$.
We will use analogous concepts for the amalgamation ${A=G_{ab}\ast_{\langle b\rangle}G_{bc}}$.
Thus we say $\mathcal{L}$ is \emph{amalgamation--reduced} if it is a sequence of triple markings alternating between
   representing elements of $G_{ab}-\langle b\rangle$ and $G_{bc}-\langle b\rangle$.
We say $\mathcal{L}$ is \emph{amalgamation--tightened} if, 
    among all amalgamation--reduced sequences $\alpha\beta\gamma$ representing the same element of $A$,
    it has maximal length prefix $\alpha$ and suffix $\gamma$ subsequences
   of triples representing elements in $\langle a\rangle\cup\langle c\rangle$.

To test $\mathcal{L}$ for niceness, 
   we modify the subroutine of 
	\cite[Lemma 20]{DLU2} to amalgamation--tighten $\mathcal{L}$.
As given, it performs Britton reduction, and then maximizes $\alpha$ length by performing a sequence of $\bigo(s)$ splittings and multiplications from left to right.  In all, it uses $\bigo(s)$ basic operations.
Afterwards, perform analogous $\langle b\rangle\langle a\rangle$ and $\langle b\rangle\langle c\rangle$ splittings and multiplications from right to left to maximize $\beta$ length. 
This does not change the $\bigo(s)$ bound on number of basic operations.
Note that $\beta\neq\emptystring$ since $\alpha\beta\gamma\notin_AC$.
Britton's Lemma implies $\mathcal{L}$ is nice unless $w$ consists of a single triple (of type $G_{ab}$ or $G_{bc}$).

There are at most $n$ intervals to test for niceness.
We can now answer whether $\widehat{x}$ or $\widehat{y}$ has a nice factor and length $\geq2$.
The original cyclic Britton--reduction, and all niceness tests, in all use $\bigo(n^2)$ basic operations.
By Lemma \ref{s steps}, we have used $\bigo(n^6)$ time so far.

Suppose, without loss of generality, that $\widehat{x}$ has length $\geq2$ and its first factor $\mathcal{L}_1$ is nice.
We want to decide whether $\widehat{x}\sim_H\widehat{y}$ in time $\bigo(n^7)$.
There are at most $n$ cyclic permutations of the intervals of $\widehat{y}$.
So, by Lemma \ref{conjugation for amalgams}, it suffices to determine in time $\bigo(n^6)$
      whether there exists $z\in C$ with $z\widehat{x}z^{-1}=\widehat{y}$.

When we tested $\mathcal{L}_1$ for niceness, we put it in the amalgamation--tightened form $\alpha\beta\gamma$.
In order to apply Lemma \ref{one candidate}, we need to compute from this a factorization $\widetilde\alpha\widetilde\beta\widetilde\gamma$
with $\widetilde\alpha,\widetilde\gamma\in C$ and $\widetilde\beta=(r_1,s_1)a^{n_1}\cdots$ tight.
If $\beta$ begins with $[u,v,w]_{ab}$, we take $\widetilde\alpha=\alpha\cdot[0,v,0]_{ab}$ and $r_1=u$.
If $\beta$ begins with $[u,v,w]_{bc}$, we take $\widetilde\alpha=\alpha$ and $r_1=2^{-v}u$.
Similarly, find $\widetilde{\alpha}'$ and $r_1'$ for the first interval $\mathcal{L}_1'$ of $\widehat{y}$.
By Lemma \ref{Lemma:reduction to tight} and Lemma \ref{one candidate}, the only candidate $z\in C$ for $z\widehat{x}=\widehat{y}z$ is
             $z=\widetilde{\alpha}'c^{r_1'-r_1}\widetilde{\alpha}^{-1}$.
				
Here $c^{r_1'-r_1}=[r_1'-r_1,0,0]_{bc}$ can be computed by the same computations underlying basic operations.
For example, if $r_1=2^{-v}u$ and $r_1'=2^{-v'}u'$, then $r_1'-r_1=2^{-(v+v')}(2^{v'}u-2^vu')$ and we can test membership in $\mathbb{Z}$.
However, we do not forget any markings, so the weight $\omega(\mathcal{T})$ increases.
Nonetheless, $\omega(\mathcal{T})=\bigo(w)$ after storing these triple markings to $\mathcal{T}$.
If $r_1'-r_1\notin\mathbb{Z}$, we immediately conclude the candidate does not work.
For notational convenience, we assume $\widetilde\alpha$ and $\widetilde\alpha'$ are empty sequences
           (by permuting the markings on $\widehat{x}$,$\widehat{y}$).

We can compute $\widehat{y}^{-1}$ from $\widehat{y}$ in time $\bigo(wn)=\bigo(n^2)$ using the fact that $[u,v,w]^{-1}=[-u,-w,-v]$.
It remains to check whether $[r_1'-r_1,0,0]_{bc}\cdot\widehat{x}\cdot[r_1-r_1',0,0]_{bc}\cdot\widehat{y}^{-1}=_H1$.
We apply again the word problem solution from \cite{DLU2}.
This sequence of triples has length $\bigo(n)$ and the weight is $\bigo(n)$, so we need an additional $\bigo(n^2)$ basic operations to solve this word problem.
The total number of basic operations up to this point is $\bigo(n^2)$, so this word problem takes time $\bigo(n^6)$.
As noted above, this means the overall algorithm takes time $\bigo(n^7)$.
\end{proof}

\begin{Rem}
In \cite[Example 3]{DMW3}, it is conjectured that the conjugacy problem can be solved in polynomial time for all \emph{hyperbolic elements} of $H=A\ast_CB$.
These are the elements not in any conjugate of $A$ or $B$.
(The name comes from the fact that they act as hyperbolic isometries on the Bass--Serre tree.)
Proposition \ref{H conjugacy runtime} leaves this conjecture unresolved.
\end{Rem}

To prove Theorem \ref{strongly generic conjugacy in H}, we will show Proposition \ref{H conjugacy runtime} solves the conjugacy problem for strongly generic input.  This requires several lemmas.

\begin{lemma}\label{Chernoff bounds}(Chernoff--Hoeffding\cite[Thm.\ 2]{Chernoff})
Suppose $X_1,X_2,\ldots$ are independent, identically distributed random variables with $-1\leq X_i\leq1$ and mean $\mu=\mu(X_i)$.
For fixed $p<\mu$:
$$\Pr\left(X_1+\cdots+X_n\leq pn\right)=e^{-\bigom(n)}.$$
\end{lemma}

\begin{lemma}\label{length not 1}(\cite[Theorem A(i) and Example 3]{DMW3})
The element of $H$ represented by a strongly generic word $w\in\Sigma^*=\{a^\pm,b^\pm,c^\pm,d^\pm\}^*$ is not in any conjugate of $A\cup B$.
In particular, $w\notin_HA\cup B$.
\end{lemma}

\begin{lemma}\label{generic with backtracking}
A strongly generic $w\in\Sigma^n=\{a^\pm,b^\pm,c^\pm,d^\pm\}^n$
   has a Britton--reduction $\widehat{w}\in(A\amalg B)^*$ with $\Omega(n)$ letters among $\{bcacb,dacad\}$.
\end{lemma}
\begin{proof}
Consider $w=t_1\cdots t_n\in\Sigma^n$.
Let $w_i=t_1\cdots t_i$.
There is a natural recursive construction of a Britton--reduction $\widehat{w_i}$ for $w_i$ so that $\widehat{w_i}$ differs from $\widehat{w_{i-1}}$ only in its rightmost letter.
For example, if $\widehat{w_{i-1}}$ ends with a letter $\sigma\in A-C$ and $t_i=d^\pm\notin A$ then $\widehat{w_i}$ is obtained by appending
$t_i$ to $\widehat{w_{i-1}}$.  On the other hand, if $\sigma\in A-C$ and $t_i\in A$ then we either replace $\sigma$ by $\sigma t_i$ (if $\sigma t_i\notin C$) or we absorb $\sigma t_i\in C$ into the previous letter of $\widehat{w}_{i-1}$.
The case $\sigma\in B-C$ is similar.
In this way, we obtain a Britton--reduction $\widehat{w}=\widehat{w_n}$.
Further, we see the \emph{Britton length} $\ell(w_i)$ of $\widehat{w_i}$ differs from $\ell(w_{i-1})$ by at most one.

We will now associate to each $w=t_1\cdots t_n\in\Sigma^n$ a string $s=s_1\cdots s_n\in\{U,D,S\}^n$.
  The letter $s_i$ will be determined by $t_i$ and the element of $H$ represented by $w_{i-1}$.
Each letter $s_i$ will have independent and identically distributed over $\Sigma^n$ probability $\Pr[U]=2/8$, $\Pr[D]=1/8$, $\Pr[S]=5/8$.

We impose requirements on the assignment $w\mapsto s_i$ when $w_{i-1}\notin_HA\cup B$.
They will be achieved in the next paragraph.
We require $s_i=U\iff\ell(w_i)>\ell(w_{i-1})$.
The connotation is that Britton length goes \emph{Up}.
Likewise, we require $s_i=D$ if $\ell(w_i)<\ell(w_{i-1})$, \emph{but not conversely}.
Thus length goes \emph{Down} or stays the same.
It follows that $\ell(w_i)=\ell(w_{i-1})$ if $s_i=S$; length stays the \emph{Same}.
If $w_{i-1}\in_HA\cup B$, we impose no requirements:
Assign the correspondence of $s_i$ to $t_i$ in an arbitrary manner achieving the desired letter distribution.

To achieve these requirements, first suppose some (hence every) Britton--reduction $\widehat{w_{i-1}}$ of $w_{i-1}$ ends in $A$.
Then $s_i=U\iff t_i=d^\pm$.  If $t_i\in\{a^\pm,c^\pm\}$ then $\ell(w_{i-1})=\ell(w_i)$ so put $s_i=S$.
At most one of $t_i=b^\pm$ can result in $\ell(w_i)<\ell(w_{i-1})$.
If it exists, this $t_i$ corresponds to $D$ and the other to $S$; otherwise, choose arbitrarily.
Similarly, if $\widehat{w_{i-1}}$ ends in $B$, take $s_i=U\iff t_i=b^\pm$ and $s_i=D$ for precisely one of $t_i=d^\pm$.

We now have the assignment $w\mapsto s_i$.
Consider strongly generic $w\in\Sigma^n$.
Lemma \ref{length not 1} implies $w_{i-1}\notin_H A\cup B$ for every $i\geq n/2$
        (simultaneously, since $n\cdot e^{-\bigom(n)}=e^{-\bigom(n)}$).
Consequently, the connotations \emph{Up}, \emph{Down or same}, and \emph{Same} apply to these $s_i$.
If $s=uv$, we say $v$ is a \emph{Dyck suffix} if each prefix of $v$ contains at least as many $U$'s as $D$'s.
We are interested in \emph{desirable} suffixes,
      namely suffixes $s_j\cdots s_n$ of $s$ with $j\geq n/2$ of the form $UUSSSSU$ followed by a Dyck suffix.
For each desirable suffix, we have $t_j\in\{b^\pm,d^\pm\}$ and the substring $t_{j+1}\cdots t_{j+5}$ gives a letter of the Britton--reduction $\widehat{w}$.
If $t_j=d^\pm$ then $t_{j+1}\cdots t_{j+5}=bcacb$ with fixed probability $p=(1/2)(1/5)^4$.
Likewise for $t_j=b^\pm$ and $t_{j+1}\cdots t_{j+5}=dacad$.
Fix some $p'\in(0,p)$.

Our strategy is to estimate the number of desirable suffixes by reading $s$ from right to left.
We have a counter $\mathcal{C}$, starting at 0, which increments (decrements) by 1 each time $U$ (resp.\ $D$) is read.
Whenever $\mathcal{C}$ ties or breaks a record maximum, the corresponding suffix of $s$ is Dyck.
We read letters of $s$ in two modes: \emph{record--breaking} (RB) mode and \emph{desirability--testing} (DT) mode.
We begin in DT mode and halt upon reading $m$ letters total in RB mode.

Whenever we enter DT mode, the suffix is Dyck ($\emptystring$ is Dyck).
Read up to 7 letters to see if we match $UUSSSSU$ (backwards), stopping early as soon as one letter mismatches.
At worst, $\mathcal{C}$ decrements by 1 each time DT mode is entered.  At best, we discover one desirable suffix.
The probability of this success is $q=(2/8)^3(5/8)^4$.  Fix $q'\in(0,q)$.
In any case, resume RB mode.

RB mode gets its own counter $\mathcal{C}'$ which behaves like $\mathcal{C}$, but only changes for letters read in RB mode.
$\mathcal{C}'$ increases on average by $1/8$ for each letter read in RB mode.
Lemma \ref{Chernoff bounds} implies that after reading $m$ letters in RB mode,
   the value of $\mathcal{C}'$ is at least $m/10$ (except with probability $e^{-\bigom(m)}$).
In particular, the $\mathcal{C}'$ record is broken at least $m/10$ times.
We enter DT mode the first $m/10$ times a $\mathcal{C}'$ record is broken.
This guarantees our precondition for entering DT mode: Each time, $\mathcal{C}$ is at a record high.

In all, we get $m/10$ chances to find a desirable suffix.
At least $mp'q'/10$ of these yield desirable suffixes contributing $bcacb$ or $dacad$ to $\widehat{w}$, by Lemma \ref{Chernoff bounds}
(except with probability $e^{-\bigom(m)}$).
We read at most $7m/10$ letters in DT mode, so at most $m+(7m/10)=17m/10$ letters overall.
Taking $m=5n/17$ ensures that we only examine suffixes with $j\geq n/2$ as required.
Thus, except with probability $e^{-\bigom(m)}+e^{-\bigom(m)}=e^{-\bigom(n)}$,
    we achieve at least $np'q'/34=\Omega(n)$ letters of $\widehat{w}$ among $\{bcacb,dacad\}$.

\end{proof}

\begin{lemma}\label{nice with backtracking}
A strongly generic word $w\in\Sigma^*$ has, after cyclic Britton--reduction, some nice factor.
\end{lemma}
\begin{proof}
Observe that $bcacb\in A$ and $dacad\in B$ from Lemma \ref{generic with backtracking} are nice.

The Britton reduction $\widehat{w}$ constructed in Lemma \ref{generic with backtracking} from 
$w=t_1\cdots t_n\in\Sigma^*$
has the property that the word $w$ decomposes as a concatenation of substrings, each of which gives a letter of $\widehat{w}$.
Every cyclic permutation of $\widehat{w}$ comes in this way from a cyclic permutation of $w$.
Therefore, by Lemma \ref{cyclic reduction amalgams}, each $w\in\Sigma^*$ has a cyclic permutation whose Britton--reduction is cyclically Britton--reduced.

By Lemma \ref{generic with backtracking}, the probability that the Britton--reduction of $w\in\Sigma^n$ has no nice factor is $e^{-\bigom(n)}$.
Since $w$ has $n$ cyclic permutations and $ne^{-\bigom(n)}=e^{-\bigom(n)}$, a strongly generic $w\in\Sigma^*$ has, after cyclic Britton--reduction, some nice factor.
\end{proof}

To relate reduced words to unreduced words, it is useful to have:

\begin{lemma}\label{Ballot Theorem}(Bertrand's Ballot Theorem\cite{Whitworth}\cite{Bertrand})
Suppose candidate $A$ gets $\alpha$ votes and $B$ gets $\beta<\alpha$ votes in an election.
The probability that $A$ always leads $B$ throughout the counting of votes is $(\alpha-\beta)/(\alpha+\beta)$.\qed
\end{lemma}

We now prove a version of Lemma \ref{nice with backtracking} for reduced words.
The proof is similar to \cite[Theorem 5]{DMW1} but in the reverse direction.

\begin{lemma}\label{nice without backtracking}
A strongly generic \emph{reduced} (or \emph{cyclically reduced}) word $w\in\Sigma^*$ has, after cyclic Britton--reduction, some nice factor.
\end{lemma}
\begin{proof}
First consider strongly generic reduced words.
Freely reducing an unreduced $w\in\Sigma^m$ involves roughly $m/8$ cancellations on average, producing a reduced word of length $3m/4$.
Therefore, to study reduced words of length $n$, we will study unreduced words of length $m=(4n/3)+i$ where we choose $i\in(-1,1]$ so
that $m$ is an integer congruent to $n$ mod 2.
To obtain a length $n$ reduced word, precisely $[m/8]$ cancellations must occur (where $[\cdot]$ is the nearest integer function).

Consider the probability $\Pr[w'\in\Sigma^n|w\in\Sigma^m]$ that exactly $[m/8]$ cancellations do occur.
This can be phrased in terms of a random walk on $\mathbb{N}$:
When at 0, we move right with probability 1.
What at $n>0$, we move left with probability $1/8$ and right with probability $7/8$.
Then $\Pr[w'\in\Sigma^n|w\in\Sigma^m]$ is the probability a random walk of length $m$ takes exactly $[m/8]$ steps to the left.
A naive estimate is:
$$\Pr[w'\in\Sigma^n|w\in\Sigma^m]\approx 
    \binom{m}{[m/8]}\left(\frac{7}8\right)^{m-[m/8]}\left(\frac{1}8\right)^{[m/8]}=\bigt\left(\frac1{\sqrt{m}}\right).$$
We use Stirling's formula $n!=\bigt((n/e)^n\cdot\sqrt{n})$.
The estimate does not account for the fact that leftward moves cannot outpace rightward moves.
Nonetheless, Bertrand's Ballot Theorem (Lemma \ref{Ballot Theorem}) gives
$$\Pr[w'\in\Sigma^n|w\in\Sigma^m]\geq\frac{m-2[m/8]}{m}\cdot 
    \binom{m}{[m/8]}\left(\frac{7}8\right)^{m-[m/8]}\left(\frac{1}8\right)^{[m/8]}=\bigt\left(\frac1{\sqrt{m}}\right).$$

Let $\mathcal{B}$ denote the set of words whose cyclic Britton--reductions have no nice factors.
Let $w'$ denote the result of freely reducing $w$.
Lemma \ref{nice with backtracking} gives:

$$\Pr[w'\in\mathcal{B}|w'\in\Sigma^n]
    \leq\frac{\Pr[w\in\mathcal{B}|w\in\Sigma^m]}{\Pr[w'\in\Sigma^n|w\in\Sigma^m]}
    \leq\frac{e^{-\Omega(m)}}{\bigom\left(\frac1{\sqrt{m}}\right)}=e^{-\Omega(n)}.$$

Note that $\Pr[\cdots|w'\in\Sigma^n]$ is the uniform distribution on reduced length $n$ words.
So we have shown a strongly generic reduced word has, after cyclic Britton--reduction, some nice factor.

The case of strongly generic \emph{cyclically reduced} words follows immediately from that of strongly generic \emph{reduced} words,
 since at least $6/7$ of all reduced words are cyclically reduced.
\end{proof}

\begin{proof}[Proof of Theorem \ref{strongly generic conjugacy in H}]
Suppose $x$ and $y$ are words in $\Sigma^*$ (respectively reduced words or cyclically reduced words) with total length $n$.
The longer, say $x$, has length $\bigt(n)$.
Then $x$ has a nice factor after cyclic Britton--reduction by Lemma \ref{nice with backtracking} (resp.\ Lemma \ref{nice without backtracking}),
   with probability $1-e^{-\Omega(n)}$.
(Further, in the strongly generic case, the cyclic Britton--reduction has length $\geq2$ by Lemma \ref{length not 1} and \cite[Theorem A(ii)]{DMW3}.)
In this case, the algorithm of Proposition \ref{H conjugacy runtime} runs in time $\bigo(n^7)$.
\end{proof}

\bibliography{higmanconjugacy}

\end{document}